\documentclass[3p,number,sort&compress]{elsarticle}

\usepackage[T1]{fontenc}
\usepackage[utf8]{inputenc}
\usepackage{lmodern}
\usepackage{microtype}
\usepackage{mathtools}
\usepackage{amssymb}
\usepackage{enumitem}
\usepackage{url}
\usepackage{balance}

\allowdisplaybreaks[3]

\newtheorem{theorem}{Theorem}
\newproof{proof}{Proof}

\begin{document}

\begin{frontmatter}

    \title{Uniqueness of solutions in thermopiezoelectricity of nonsimple materials}


 	\author[2]{Francesca Passarella}
 	\ead{fpassarella@unisa.it}

 	\author[2]{Vincenzo Tibullo\corref{cor1}}
 	\ead{vtibullo@unisa.it}

	\cortext[cor1]{Corresponding author}

	\affiliation[2]{organization={Dipartimento di Matematica, Università di Salerno}}

	\begin{abstract}
		In this paper, using the entropy production inequality of Green and Laws, the theory of thermopiezoelectricity is presented, in which the second gradient of displacement and the second gradient of electric potential are included in the set of independent constitutive variables.
		At first, the appropriate thermodynamic restrictions and constitutive equations are obtained.
		Then, the basic equations of linear theory are derived and a uniqueness result for the mixed boundary-initial value problem is presented.
	\end{abstract}

	
	\date{March 2, 2022}

\end{frontmatter}

\section{Introduction}

The classical theory of heat propagation is based on two equations: the energy balance equation and the Fourier law for the heat flux.
These equations lead to the classical heat equation that has the unrealistic feature that the velocity of heat propagation is infinite.
This classical theory was firstly questioned by \citet{Cattaneo1948, Cattaneo1958} and \citet{Vernotte1958}, leading to the rise of new fields of research: heat waves and second sound propagation can be found in the work of \citet{Straughan2011b} and in extended irreversible thermodynamics (EIT) in the books of \citet{Jou1996} and \citet{Mueller1998}.
In~\cite{Gurtin1968} \citeauthor{Gurtin1968} considered a general theory of heat conduction with finite velocity, based on a heat flux law with memory.
For more information on thermoelasticity theories that predict a finite velocity for the propagation of thermal signals, see the the reviews of \citet{Chandrasekharaiah1998, Chandrasekharaiah1986}.
Furthermore, as presented by \citeauthor{Iesan2004c} in~\cite{Iesan2004c}, the theory proposed by \citet{Green1972a} is an alternative way of formulation of the propagation of heat.
They make use of an entropy inequality in which a new constitutive function appears (see e.g.~\cite{Passarella2013}).

The origin of the theory of nonsimple elastic materials goes back to the works of \citet{Toupin1962, Toupin1964} and \citet{Mindlin1964}.
The first investigations of nonsimple thermoelastic materials are presented in~\cite{Ahmadi1975, Batra1976}.
Within the framework of that theory for this materials, the second-order displacement gradient is added to the classical set of independent constitutive variables.
The theory of nonsimple thermoelastic materials has been discussed in various papers (see for example~\cite{Mindlin1968, Ahmadi1977, Iesan1983, Ciarletta1989a, Kalpakides2002, Aouadi2019, Aouadi2020, Aouadi2020a, Bartilomo1997})

The problem of the interaction of electromagnetic fields with elastic solids was the subject of important investigations (see e.g.~\cite{Eringen2004, Truesdell1960, Parkus1972, Grot1976, Nowacki1983, Maugin1988} and the literature cited therein).
Certain crystals (for example quartz) when subject to stress, become electrically polarized (piezoelectric effect).
Conversely, an external electromagnetic field produces deformation in a piezoelectric crystal.
The theory of thermopiezoelectricity has been studied in various works (see e.g.~\cite{Mindlin1961, Eringen2012, Iesan2018, Nowacki1978, Morro1991}).

In this paper, we derive a theory of thermopiezoelectricity of a body in which the second gradient of displacement and the second gradient of electric potential are included in the set of independent constitutive variables.
We obtain the appropriate thermodynamic restrictions and constitutive equations, with the help of an entropy production inequality proposed by \citet{Green1972a}.
Next, for both anisotropic and isotropic materials we establish the basic equations of the linear theory and we obtain a uniqueness result for the mixed initial-boundary values problem.

\section{Basic equations}

We consider a body that at some instant occupies the region \(B\) of the Euclidean three-dimensional space and is bounded by the piecewise smooth surface \(\partial B\).
The motion of the body is referred to the reference configuration \(B\) and to a fixed system of rectangular Cartesian axes \(Ox_i\) \((i=1,2,3)\).

We shall employ the usual summation and differentiation conventions: Latin subscripts are understood to range over the integers \((1,2,3)\), summation over repeated subscripts is implied and subscripts preceded by a comma denote partial differentiation with respect to the corresponding Cartesian coordinate.
In what follows we use a superposed dot to denote partial differentiation with respect to the time \(t\).
Further, we will neglect the issues of regularity, simply understanding a degree of smoothness sufficient to make sense everywhere.

We restrict our attention to the theory of homogeneous piezoelectric solids.
Let \(u_i\) be the displacement vector, \(D_i\) the electric displacement vector and \(E_i\) the electric field vector.
The equation of motion is
\begin{equation}
	t_{ji,j} + \rho f_i = \rho \ddot{u}_i,
	\label{eq:motion}
\end{equation}
while the equations for the quasi-static electric field are~\cite{Eringen2012}
\begin{equation}
	D_{j,j} = g, \qquad E_i = -\varphi_{,i},
	\label{eq:electric}
\end{equation}
where \(t_{ij}\) is the stress tensor, \(f_i\) is the external body force per unit mass, \(\rho\) is the reference mass density, \(g\) is the density of free charges and \(\varphi\) is the electric potential.

Following~\cite{Toupin1964, Eringen2004, Iesan2018}, we assume that the body is free from initial stress and we postulate an energy balance in the form
\begin{alignat*}{1}
	 & \frac{d}{dt} \int_P \rho \left( e + \frac{1}{2}\dot{u}_i\dot{u}_i \right) dv = \int_P ( \rho f_i\dot{u}_i + \dot{D}_i E_i + \rho r ) dv \\
	 & \qquad + \int_{\partial P} ( t_i\dot{u}_i + \mu_{ji}\dot{u}_{i,j} + \dot{Q}_i E_i-q ) da
\end{alignat*}
for every part \(P\) of \(B\) and for every time.
Here \(e\) is the internal energy per unit mass, \(r\) is the heat supply per unit mass.
Moreover, \(t_i\) is the traction vector, \(\mu_{ji}\) is the hypertraction tensor, \(Q_i\) is the generalized surface charge density and \(q\) is the heat flux.

Following \citet{Iesan2004c, Iesan2018}, for the traction vector, the hyperstress tensor \(\mu_{kji}\), the heat flux vector \(q_i\) and the electric quadrupole \(Q_{ji}\) we have
\begin{equation}
	t_i = t_{ji} n_j, \quad \mu_{ji} = \mu_{kji}n_k, \quad q = q_i n_i, \quad \dot{Q}_i = \dot{Q}_{ji}n_j,
	\label{eq:tensors}
\end{equation}
where \(n_i\) is the outward unit normal vector to the boundary surface \(\partial B\).

By using eqs.~\eqref{eq:motion} and~\eqref{eq:tensors} and the divergence theorem, the balance of energy becomes
\begin{alignat*}{1}
	 \int_P \rho\dot{e} dv &= \int_P [ ( t_{ji} + \mu_{kji,k} ) \dot{u}_{i,j} + \mu_{kji}\dot{u}_{i,jk} \\
	 & + ( \dot{D}_i + \dot{Q}_{ji,j} ) E_i + \dot{Q}_{ji}E_{i,j} + \rho r - q_{i,i} ] dv.
\end{alignat*}
If we introduce the tensors
\begin{equation}
	\tau_{ji} = t_{ji} + \mu_{kji,k}, \qquad \sigma_i = D_i + Q_{ji,j},
	\label{eq:tau}
\end{equation}
eqs.~\eqref{eq:motion} and~\eqref{eq:electric}\(_1\) can be written as
\begin{equation}
	\tau_{ji,j} - \mu_{kji,kj} + \rho f_i = \rho\ddot{u}_i, \qquad \sigma_{i,i} - Q_{ji,ji} = g,
	\label{eq:motion-final}
\end{equation}
and the local form of energy balance becomes
\begin{equation}
	\rho\dot{e} = \tau_{ji}\dot{u}_{i,j} + \mu_{kji}\dot{u}_{i,jk} + \dot{\sigma}_i E_i + \dot{Q}_{ji}E_{i,j} + \rho r - q_{i,i}.
	\label{eq:energy-local}
\end{equation}
As proved by \citet{Iesan2004c}, the invariance of the energy equation for observers in rotating motion the one with respect to the other, leads to
\[
	\tau_{ji} = \tau_{ij},
\]
moreover, without loss of generality, we can suppose
\[
	\mu_{kji} = \mu_{jki},
\]
given that the skew symmetric part makes no contribution to the rate of work over any closed surface in the body, or over the boundary.

Let's set
\begin{equation}
	e_{ij} = \frac{1}{2} ( u_{i,j} + u_{j,i} ) , \quad \kappa_{ijk} = u_{k,ij}, \quad V_{ij} = E_{i,j},
	\label{eq:geometric}
\end{equation}
then eq.~\eqref{eq:energy-local} becomes
\begin{equation}
	\rho\dot{e} = \tau_{ji}\dot{e}_{ij} + \mu_{ijk}\dot{\kappa}_{ijk} + \dot{\sigma}_i E_i + \dot{Q}_{ji}V_{ij} + \rho r - q_{i,i}.
	\label{eq:dot-energy}
\end{equation}

We postulate the entropy production inequality proposed by \citet{Green1972a}
\[
	\frac{d}{dt} \int_P \rho\eta dv \geq \int_P \frac{\rho r}{\phi} dv - \int_{\partial P} \frac{q}{\phi} da
\]
for every part \(P\) of \(B\) and every time.
Here \(\eta\) is he entropy per unit mass, \(\phi\) is a new strictly positive thermal variable needing a constitutive equation.

Using eq.~\eqref{eq:tensors} and the divergence theorem, we have in local form
\[
	\rho\dot{\eta} \geq \frac{\rho r}{\phi} - \left( \frac{q_i}{\phi} \right)_{,i}
\]
or equivalently
\begin{equation}
	\rho\phi\dot{\eta} \geq \rho r - q_{i,i} + \frac{1}{\phi}q_i\phi_{,i}.
	\label{eq:entropy}
\end{equation}

We now introduce the scalar function \(\sigma\) (electric enthalpy) defined by
\begin{equation}
	\sigma = \rho ( e - \phi\eta ) - \sigma_i E_i - Q_{ji}V_{ij},
	\label{eq:enthalpy}
\end{equation}
and we substitute eqs.~\eqref{eq:dot-energy} and~\eqref{eq:enthalpy} in the inequality~\eqref{eq:entropy}, we then obtain
\begin{equation}
	\dot{\sigma} + \rho\dot{\phi}\eta - \tau_{ji}\dot{e}_{ij} - \mu_{ijk}\dot{\kappa}_{ijk} + \sigma_i\dot{E}_i + Q_{ji}\dot{V}_{ij} + \frac{1}{\phi}q_i\phi_{,i} \leq 0.
	\label{eq:a-entropy}
\end{equation}

\section{Thermodynamics restrictions}

Let \(\theta\) be the difference of the absolute temperature \(T\) and the absolute temperature in the reference configuration \(T_0\), i.e. \(\theta = T-T_0\).

We require constitutive equations for \(\sigma\), \(\phi\), \(\eta\), \(\tau_{ij}\), \(\mu_{ijk}\), \(\sigma_i\), \(Q_{ij}\) and \(q_i\) and assume that these are functions of the set of variables \(\mathcal{S} = \left(e_{ij},\kappa_{ijk},E_i,\theta,\dot{\theta},\theta_{,i},V_{ij}\right)\)
\begin{equation}
	\begin{alignedat}{3}
		& \sigma = \sigma(\mathcal{S}), \quad && \phi = \phi(\mathcal{S}), \quad && \eta = \eta(\mathcal{S}), \\
		& \tau_{ij} = \tau_{ij}(\mathcal{S}), \quad && \mu_{ijk} = \mu_{ijk}(\mathcal{S}), \quad && \sigma_i = \sigma_i(\mathcal{S}), \\
		& Q_{ij} = Q_{ij}(\mathcal{S}), \quad && q_i = q_i(\mathcal{S}),
	\end{alignedat}
	\label{eq:constitutive-eqs}
\end{equation}
where \(e_{ij}=e_{ji}\), \(\kappa_{ijk}=\kappa_{jik}\) e \(V_{ij}=V_{ji}\).

If we replace eqs.~\eqref{eq:constitutive-eqs} into the inequality~\eqref{eq:a-entropy}, we arrive to
{\small
	\begin{alignat*}{1}
		 & \left[ \frac{\partial\sigma}{\partial\theta} + \rho\eta \frac{\partial\phi}{\partial\theta} \right] \dot{\theta} + \left[ \frac{\partial\sigma}{\partial\dot{\theta}} + \rho\eta \frac{\partial\phi}{\partial\dot{\theta}} \right] \ddot{\theta} \\
		 & + \left[ \frac{\partial\sigma}{\partial e_{ij}} + \rho\eta \frac{\partial\phi}{\partial e_{ij}} - \tau_{ji} \right] \dot{e}_{ij} + \left[ \frac{\partial\sigma}{\partial\kappa_{ijk}} + \rho\eta \frac{\partial\phi}{\partial\kappa_{ijk}} - \mu_{ijk} \right] \dot{\kappa}_{ijk} \\
		 & + \left[ \frac{\partial\sigma}{\partial\theta_{,i}} + \rho\eta \frac{\partial\phi}{\partial\theta_{,i}} + \frac{1}{\phi} q_i \frac{\partial\phi}{\partial\dot{\theta}} \right] \dot{\theta}_{,i} + \left[ \frac{\partial\sigma}{\partial E_i} + \rho\eta \frac{\partial\phi}{\partial E_i} + \sigma_i \right] \dot{E_i} \\
		 & + \left[ \frac{\partial\sigma}{\partial V_{ij}} + \rho\eta \frac{\partial\phi}{\partial V_{ij}} + Q_{ji} \right] \dot{V}_{ij} + \frac{1}{\phi} q_l \bigg[ \frac{\partial\phi}{\partial e_{ij}} e_{ij,l} + \frac{\partial\phi}{\partial\kappa_{ijk}} \kappa_{ijk,l} \\
		 & + \frac{\partial\phi}{\partial\theta} \theta_{,l} + \frac{\partial\phi}{\partial\theta_{,i}} \theta_{,il} + \frac{\partial\phi}{\partial E_i} E_{i,l} + \frac{\partial\phi}{\partial V_{ij}} V_{ij,l} \bigg] \leq 0.
	\end{alignat*}
}%
We have the following conditions
\begin{alignat*}{2}
	q_l \frac{\partial\phi}{\partial\kappa_{ijk}} &= 0, \quad & q_l \frac{\partial\phi}{\partial\theta_{,i}} + q_i \frac{\partial\phi}{\partial\theta_{,l}} &= 0, \\
	q_l \frac{\partial\phi}{\partial V_{ij}} &= 0 \quad & q_l \frac{\partial\phi}{\partial E_i} + q_i \frac{\partial\phi}{\partial E_l} &= 0,
\end{alignat*}
and if we suppose that \(q_i \neq 0\), we obtain
\[
	\frac{\partial\phi}{\partial\kappa_{ijk}} = 0, \quad \frac{\partial\phi}{\partial\theta_{,i}} = 0, \quad \frac{\partial\phi}{\partial E_i} = 0, \quad \frac{\partial\phi}{\partial V_{ij}} = 0,
\]
so that \(\phi\) can only depend on \(e_{ij}\), \(\theta\), \(\dot{\theta}\), i.e.
\begin{equation}
	\phi = \phi ( e_{ij}, \theta, \dot{\theta} ).
	\label{eq:phi-dependence}
\end{equation}
Consequently, assuming \(\partial\phi/\partial\dot{\theta} \neq 0\), we have the following  restrictions
\begin{equation}
	\begin{alignedat}{1}
		\tau_{ij} &= \frac{\partial\sigma}{\partial e_{ij}} + \rho\eta \frac{\partial\phi}{\partial e_{ij}} \\
		\mu_{ijk} &= \frac{\partial\sigma}{\partial\kappa_{ijk}} \\
		\rho\eta &= -\left( \frac{\partial\phi}{\partial\dot{\theta}} \right)^{-1} \frac{\partial\sigma}{\partial\dot{\theta}} \\
		q_i &= -\phi \left( \frac{\partial\phi}{\partial\dot{\theta}} \right)^{-1} \frac{\partial\sigma}{\partial\theta_{,i}} \\
		\sigma_i &= -\frac{\partial\sigma}{\partial E_i}, \\
		Q_{ji} &= -\frac{\partial\sigma}{\partial V_{ij}}
	\end{alignedat}
	\label{eq:stress-definitions}
\end{equation}
and the dissipation inequality is
\begin{equation}
	\left( \frac{\partial\sigma}{\partial\theta} + \rho\eta \frac{\partial\phi}{\partial\theta} \right) \dot{\theta} + \frac{1}{\phi} q_l \left( \frac{\partial\phi}{\partial\theta} \theta_{,l} + \frac{\partial\phi}{\partial e_{ij}} \kappa_{lij} \right) \leq 0.
	\label{eq:entropy-reduced}
\end{equation}

Taking into account eqs.~\eqref{eq:enthalpy},~\eqref{eq:phi-dependence},~\eqref{eq:stress-definitions}, the equation of energy~\eqref{eq:dot-energy} reduces to
\[
	\phi\rho\dot{\eta} + q_{i,i} - \rho r + \left( \frac{\partial\sigma}{\partial\theta} + \rho\eta \frac{\partial\phi}{\partial\theta} \right) \dot{\theta} + \frac{\partial\sigma}{\partial\theta_{,i}} \dot{\theta}_{,i} = 0.
\]

\section{Linear constitutive equations}

We consider a quadratic Taylor expansion for \(\sigma\) and \(\phi\) with initial point the reference state
{\small
	\begin{equation}
		\begin{alignedat}{1}
			 \sigma &= \alpha_{ij}^{(1)}e_{ij} + \alpha_{ijk}^{(2)}\kappa_{ijk} + \alpha_i^{(3)}E_i + \alpha^{(4)}\theta + \alpha^{(5)}\dot{\theta} + \alpha_i^{(6)}\theta_{,i} \\
			 & + \alpha_{ij}^{(7)}V_{ij} + \frac{1}{2} \Bigl( a_{ijkl}^{(11)}e_{ij}e_{kl} + a_{ijklhm}^{(22)}\kappa_{ijk}\kappa_{lhm} + a_{ij}^{(33)}E_i E_j \\
			 & + a^{(44)}\theta^2 + a^{(55)}\dot{\theta}^2 + a_{ij}^{(66)}\theta_{,i}\theta_{,j} + a_{ijkh}^{(77)}V_{ij}V_{hk} \Bigr) \\
			 & + a_{ijklh}^{(12)}e_{ij}\kappa_{klh} + a_{ijk}^{(13)}e_{ij}E_k + a_{ij}^{(14)}e_{ij}\theta + a_{ij}^{(15)}e_{ij}\dot{\theta} \\
			 & + a_{ijk}^{(16)}e_{ij}\theta_{,k} + a_{ijkl}^{(17)}e_{ij}V_{kl} + a_{ijkl}^{(23)}\kappa_{ijk}E_l + a_{ijk}^{(24)}\kappa_{ijk}\theta \\
			 & + a_{ijk}^{(25)}\kappa_{ijk}\dot{\theta} + a_{ijkl}^{(26)}\kappa_{ijk}\theta_{,l} + a_{ijkhl}^{(27)}\kappa_{ijk}V_{hl} + a_i^{(34)}E_i\theta \\
			 & + a_i^{(35)}E_i\dot{\theta} + a_{ij}^{(36)}E_i\theta_{,j} + a_{ijk}^{(37)}E_i V_{jk} + a^{(45)}\theta\dot{\theta} + a_i^{(46)}\theta\theta_{,i} \\
			 & + a_{ij}^{(47)}\theta V_{ij} + a_i^{(56)}\dot{\theta}\theta_{,i} + a_{ij}^{(57)}\dot{\theta}V_{ij} + a_{ijk}^{(67)}\theta_{,i}V_{jk}
		\end{alignedat}
		\label{eq:sigma-quadratic-taylor}
	\end{equation}
}%
and
\begin{equation}
	\begin{alignedat}{1}
		\phi &= T_0 + \beta_{ij}^{(1)}e_{ij} + \beta^{(4)}\theta + \beta^{(5)}\dot{\theta} \\
		& + \frac{1}{2} \left( b_{ijkl}^{(11)}e_{ij}e_{kl} + b^{(44)}\theta^2 + b^{(55)}\dot{\theta}^2 \right) \\
		& + b_{ij}^{(14)}e_{ij}\theta + b_{ij}^{(15)}e_{ij}\dot{\theta} + b^{(45)}\theta\dot{\theta}.
	\end{alignedat}
	\label{eq:phi-quadratic-taylor}
\end{equation}

The coefficients in~\eqref{eq:sigma-quadratic-taylor} and~\eqref{eq:phi-quadratic-taylor} satisfy the following symmetry relations
\begin{equation}
	\begin{alignedat}{2}
		& \alpha_{ij}^{(1)} = \alpha_{ji}^{(1)}, \quad && \alpha_{ijk}^{(2)} = \alpha_{jik}^{(2)}, \\
		& \alpha_{ij}^{(7)} = \alpha_{ji}^{(7)} \quad && a_{ijkl}^{(11)} = a_{jikl}^{(11)} = a_{klij}^{(11)}, \\
		& a_{ijklhm}^{(22)} = a_{jiklhm}^{(22)} = a_{lhmijk}^{(22)}, \quad && a_{ij}^{(33)} = a_{ji}^{(33)}, \\
		& a_{ij}^{(66)} = a_{ji}^{(66)}, \quad && a_{ijkl}^{(77)} = a_{jikl}^{(77)} = a_{klij}^{(77)}, \\
		& a_{ijklh}^{(12)} = a_{jiklh}^{(12)} = a_{ijlkh}^{(12)}, \quad && a_{ijk}^{(13)} = a_{jik}^{(13)}, \\
		& a_{ij}^{(14)} = a_{ji}^{(14)}, \quad && a_{ij}^{(15)} = a_{ji}^{(15)}, \\
		& a_{ijk}^{(16)} = a_{jik}^{(16)}, \quad && a_{ijkl}^{(17)} = a_{jikl}^{(17)} = a_{klij}^{(17)}, \\
		& a_{ijkl}^{(23)} = a_{jikl}^{(23)}, \quad && a_{ijk}^{(24)} = a_{jik}^{(24)}, \\
		& a_{ijk}^{(25)} = a_{jik}^{(25)}, \quad && a_{ijkl}^{(26)} = a_{jikl}^{(26)}, \\
		& a_{ijkhl}^{(27)} = a_{jikhl}^{(27)} = a_{ijklh}^{(27)}, \quad && a_{ijk}^{(37)} = a_{ikj}^{(37)}, \\
		& a_{ij}^{(47)} = a_{ji}^{(47)}, \quad && a_{ij}^{(57)} = a_{ji}^{(57)}, \\
		& a_{ijk}^{(67)} = a_{ikj}^{(67)}, \\
		& \beta_{ij}^{(1)} = \beta_{ji}^{(1)}, \quad && b_{ijkl}^{(11)} = b_{jikl}^{(11)} = b_{klij}^{(11)}, \\
		& b_{ij}^{(14)} = b_{ji}^{(14)}, \quad && b_{ij}^{(15)} = b_{ji}^{(15)}.
	\end{alignedat}
	\label{eq:symmetries}
\end{equation}

We suppose that the positive function \(\phi(e_{ij},\theta,\dot{\theta})\) satisfies
\[
	\phi = \phi ( e_{ij}, \theta, 0) = T_0 + \theta
\]
and, consequently, we have
\begin{alignat*}{3}
	 & \beta_{ij}^{(1)} = 0, \quad && \beta^{(4)} = 1, \\
	 & b_{ijkl}^{(11)} = 0, \quad && b^{(44)} = 0, \quad && b_{ij}^{(14)} = 0,
\end{alignat*}
so that
\[
	\phi = T_0 + \theta + \beta\dot{\theta} + \frac{1}{2} b^{(55)}\dot{\theta}^2 + b_{ij}^{(15)}e_{ij}\dot{\theta} + b^{(45)}\theta\dot{\theta}.
\]

Assuming that the body in the reference state is free from stress and hyperstress and has zero electric displacement vector, heat flux and electric quadrupole, we obtain
\begin{alignat*}{3}
	 & \alpha_{ij}^{(1)} = 0, \quad && \alpha_{ijk}^{(2)} = 0, \\
	 & \alpha_i^{(3)} = 0, \quad && \alpha_i^{(6)} = 0, \quad && \alpha_{ij}^{(7)} = 0,
\end{alignat*}
moreover, linearizing eqs.~\eqref{eq:stress-definitions} we arrive to
\begin{equation}
	\begin{alignedat}{1}
		 \rho\eta &= -\frac{1}{\beta} \Bigl( \alpha^{(5)} + a_{ij}^{(15)}e_{ij} + a_{ijk}^{(25)}\kappa_{ijk} + a_i^{(35)}E_i \\
		 & + a^{(45)}\theta + a^{(55)}\dot{\theta} + a_i^{(56)}\theta_{,i} + a_{ij}^{(57)}V_{ij} \Bigr) \\
		 & + \frac{\alpha^{(5)}}{\beta^2} \left( b_{ij}^{(15)}e_{ij} + b^{(45)}\theta + b^{(55)}\dot{\theta} \right),
	\end{alignedat}
	\label{eq:rho-eta}
\end{equation}
and
{\small
	\begin{equation}
		\begin{alignedat}{1}
			\tau_{ij} &= a_{ijkl}^{(11)}e_{kl} + a_{ijklh}^{(12)}\kappa_{klh} + a_{ijk}^{(13)}E_k + a_{ij}^{(14)}\theta \\
			& +\Big( a_{ij}^{(15)} - \frac{\alpha^{(5)}}{\beta} b_{ij}^{(15)} \Big) \dot{\theta} + a_{ijk}^{(16)}\theta_{,k} + a_{ijkl}^{(17)}V_{kl}, \\
			\mu_{ijk} &= a_{hlijk}^{(12)}e_{hl} + a_{ijkhlm}^{(22)}\kappa_{hlm} + a_{ijkh}^{(23)}E_h + a_{ijk}^{(24)}\theta \\
			& +a_{ijk}^{(25)}\dot{\theta} + a_{ijkh}^{(26)}\theta_{,h} + a_{ijklh}^{(27)}V_{lh}, \\
			\sigma_i &= -\big( a_{jki}^{(13)}e_{jk} + a_{jkli}^{(23)}\kappa_{jkl} + a_{ij}^{(33)}E_j + a_i^{(34)}\theta \\
			& +a_i^{(35)}\dot{\theta} + a_{ij}^{(36)}\theta_{,j} + a_{ijk}^{(37)}V_{jk} \big), \\
			q_i &= -\frac{T_0}{\beta} \big( a_{jki}^{(16)}e_{jk} + a_{jkli}^{(26)}\kappa_{jkl} + a_{ji}^{(36)}E_j + a_i^{(46)}\theta \\
			& +a_i^{(56)}\dot{\theta} + a_{ij}^{(66)}\theta_{,j} + a_{ijk}^{(67)}V_{jk} \big), \\
			Q_{ij} &= -\big( a_{klij}^{(17)}e_{kl} + a_{klmij}^{(27)}\kappa_{klm} + a_{kij}^{(37)}E_k + a_{ij}^{(47)}\theta \\
			& +a_{ij}^{(57)}\dot{\theta} + a_{kij}^{(67)}\theta_{,k} + a_{ijkl}^{(77)}V_{kl} \big).
		\end{alignedat}
		\label{eq:constitutive}
	\end{equation}
}%

Linearizing the coefficients of the entropy inequality~\eqref{eq:entropy-reduced} we obtain
\[
	\mathcal{A}(\mathcal{S})\dot{\theta} + \mathcal{B}_i(\mathcal{S})\theta_{,i} \leq 0
\]
where
\begin{alignat*}{1}
	\mathcal{A}(\mathcal{S}) &= \alpha^{(4)} + a_{ij}^{(14)}e_{ij} + a_{ijk}^{(24)}\kappa_{ijk} + a_i^{(34)}E_i \\
    & +a^{(44)}\theta + a^{(45)}\dot{\theta} + a_i^{(46)}\theta_{,i} + a_{ij}^{(47)}V_{ij} \\
    & -\frac{1}{\beta} \big( \alpha^{(5)} + a_{ij}^{(15)}e_{ij} + a_{ijk}^{(25)}\kappa_{ijk} + a_i^{(35)}E_i \\
    & +a^{(45)}\theta + a^{(55)}\dot{\theta} + a_i^{(56)}\theta_{,i} + a_{ij}^{(57)}V_{ij} \big) \\
    & +\frac{\alpha^{(5)}}{\beta^2} \big( b_{ij}^{(15)}e_{ij} + b^{(45)} ( \theta - \beta\dot{\theta} ) + b^{(55)}\dot{\theta} \big), \\
	\mathcal{B}_i(\mathcal{S}) &= -\frac{1}{\beta} \big( a_{jki}^{(16)}e_{jk} + a_{jkli}^{(26)}\kappa_{jkl} + a_{ji}^{(36)}E_j + a_i^{(46)}\theta \\
    & +a_i^{(56)}\dot{\theta} + a_{ij}^{(66)}\theta_{,j} + a_{ijk}^{(67)}V_{jk} \big).
\end{alignat*}

If we define
\[
	\mathcal{S}_0 = \mathcal{S} \big|_{\dot{\theta} = 0, \theta_{,i} = 0} = \left( e_{ij}, \kappa_{ijk}, E_i, \theta, 0, 0, V_{ij} \right)
\]
by virtue of the dissipation inequality we arrive to
\[
	\mathcal{A}(\mathcal{S}_0) = 0, \qquad \mathcal{B}_i(\mathcal{S}_0) = 0.
\]
From the first condition we have
\[
	\alpha^{(5)} = \beta\alpha^{(4)}
\]
and
\begin{alignat*}{2}
	a_{ij}^{(15)} &=\beta a_{ij}^{(14)} + \alpha^{(4)}b_{ij}^{(15)}, \\
	a_{ijk}^{(25)} &= \beta a_{ijk}^{(24)}, \quad & a_i^{(35)} &= \beta a_i^{(34)}, \quad \\
	a^{(45)} &= \beta a^{(44)} + \alpha^{(4)}b^{(45)}, \quad & a_{ij}^{(57)} &= \beta a_{ij}^{(47)},
\end{alignat*}
while from the second one we obtain
\begin{alignat*}{3}
	a_{ijk}^{(16)} &= 0, \quad & a_{ijkl}^{(26)} &= 0, \quad & a_{ij}^{(36)} &= 0, \\
	a_i^{(46)} &= 0, \quad & a_{ijk}^{(67)} &= 0.
\end{alignat*}

Finally we have
\begin{alignat*}{1}
	\mathcal{A}(\mathcal{S}) &= -\frac{1}{\beta} ( \gamma\dot{\theta} + a_i^{(56)}\theta_{,i} ), \\
	\mathcal{B}_i(\mathcal{S}) &= -\frac{1}{\beta} ( a_i^{(56)}\dot{\theta} + a_{ij}^{(66)}\theta_{,j} ),
\end{alignat*}
where we defined
\[
	\gamma = a^{(55)} - \beta^2 a^{(44)} - \alpha^{(4)}b^{(55)}.
\]

Using these relations, the constitutive equations~\eqref{eq:rho-eta} and~\eqref{eq:constitutive} can be expressed as
\begin{equation}
	\begin{alignedat}{1}
		\tau_{ij} &= a_{ijkl}^{(11)}e_{kl} + a_{ijklh}^{(12)}\kappa_{klh} + a_{ijk}^{(13)}E_k \\
	    & +a_{ij}^{(14)} ( \theta + \beta\dot{\theta} ) + a_{ijkl}^{(17)}V_{kl}, \\
		\mu_{ijk} &= a_{lhijk}^{(12)}e_{lh} + a_{ijklhm}^{(22)}\kappa_{lhm} + a_{ijkl}^{(23)}E_l \\
        & +a_{ijk}^{(24)} ( \theta + \beta\dot{\theta} ) + a_{ijklh}^{(27)}V_{lh}, \\
		-\sigma_i &= a_{jki}^{(13)}e_{jk} + a_{jkli}^{(23)}\kappa_{jkl} + a_{ij}^{(33)}E_j \\
        & +a_i^{(34)} ( \theta + \beta\dot{\theta} ) + a_{ijk}^{(37)}V_{jk}, \\
		-\rho\eta &= a_{ij}^{(14)}e_{ij} + a_{ijk}^{(24)}\kappa_{ijk} + a_i^{(34)}E_i \\
        & +a^{(44)} ( \theta + \beta\dot{\theta} ) + a_{ij}^{(47)}V_{ij} \\
        & +\frac{1}{\beta} \left( \gamma\dot{\theta} + a_i^{(56)}\theta_{,i} \right) + \alpha^{(4)}, \\
		-Q_{ij} &= a_{klij}^{(17)}e_{kl} + a_{klmij}^{(27)}\kappa_{klm} + a_{kij}^{(37)}E_k \\
        & +a_{ij}^{(47)} ( \theta + \beta\dot{\theta} ) + a_{ijkl}^{(77)}V_{kl}, \\
		-\frac{q_i}{T_0} &= \frac{1}{\beta} \left( a_i^{(56)}\dot{\theta} + a_{ij}^{(66)}\theta_{,j} \right).
	\end{alignedat}
	\label{eq:cost-eq}
\end{equation}

In the linear approximation, the energy equation is
\begin{equation}
	\rho T_0\dot{\eta} + q_{i,i} - \rho r = 0
	\label{eq:dot-energy-linear}
\end{equation}
or, with help of eqs.~\eqref{eq:cost-eq}, we arrive to
\begin{align*}
	 & a_{ij}^{(14)}\dot{e}_{ij} + a_{ijk}^{(24)}\dot{\kappa}_{ijk} + a_i^{(34)}\dot{E}_i + a^{(44)} ( \dot{\theta} + \beta\ddot{\theta} ) + a_{ij}^{(47)}\dot{V}_{ij} \\
	 & \qquad + \frac{1}{\beta} \left( \gamma\ddot{\theta} +2a_i^{(56)}\dot{\theta}_{,i} + a_{ij}^{(66)}\theta_{,ij} \right) + \frac{1}{T_0} \rho r = 0.
\end{align*}
The dissipation inequality~\eqref{eq:entropy-reduced} becomes
\[
	\frac{1}{\beta} \left( \gamma\dot{\theta}^2 + 2a_i^{(56)}\dot{\theta}\theta_{,i} + a_{ij}^{(66)}\theta_{,i}\theta_{,j} \right) \geq 0.
\]
The following quadratic form is defined from the previous dissipation inequality
\begin{equation}
	\mathcal{P} ( \xi, \eta_i ) = \frac{1}{\beta} \left( \gamma\xi^2 + 2a_i^{(56)}\xi\eta_i + a_{ij}^{(66)}\eta_i\eta_j \right)
	\label{eq:quadratic-form-p}
\end{equation}
so that
\begin{equation}
	\mathcal{P} ( \xi, \eta_i ) \geq 0, \qquad \forall \xi, \eta_i.
	\label{eq:semidef}
\end{equation}
If we consider the symmetric matrix associated to the quadratic form \(\mathcal{P}\)
\[
	\frac{1}{\beta}
	\begin{pmatrix}
		\gamma & a_j^{(56)} \\
		a_i^{(56)} & a_{ij}^{(66)}
	\end{pmatrix},
\]
its positive semi-definiteness implies, in particular, that
\[
	\frac{\gamma}{\beta} \geq 0, \qquad \frac{1}{\beta} a_{ij}^{(66)}\eta_i\eta_j \geq 0, \qquad \forall \eta_i.
\]

The basic equations of linear theory of thermopiezoelectric solids consist of the equations of motion~\eqref{eq:motion-final}, the equation of energy~\eqref{eq:dot-energy-linear}, the geometrical equations~\eqref{eq:geometric},~\eqref{eq:electric}\(_2\), the constitutive equations~\eqref{eq:cost-eq} with the restriction~\eqref{eq:semidef}, on \(B\times I\), where \(I=[0,t_0)\), where \(t_0\leq+\infty\). 

Following \citet{Toupin1964} and \citet{Mindlin1964}, we consider \(P_i\), \(R_i\), \(\Lambda\) and \(H\) defined in such a way that the total rate of work of the surface forces over the smooth surface \(\partial P\) can be expressed in the form
\begin{alignat*}{1}
	 & \int_{\partial P} \left( t_{ki}\dot{u}_i + \mu_{kji}\dot{u}_{i,j} - \varphi\dot{D}_k - \varphi_{,i}\dot{Q}_{ki} \right) n_k dA \\
	 & \qquad = \int_{\partial P} \left( P_i\dot{u}_i + R_i\mathcal{D}\dot{u}_i - \varphi\dot{\Lambda} - \dot{H}\mathcal{D}\varphi \right) dA.
\end{alignat*}
Here we used
\begin{alignat*}{1}
	& P_i = ( \tau_{ji} - \mu_{kji,k} ) n_j - \mathcal{D}_j ( \mu_{kji} n_k ) + ( \mathcal{D}_l n_l ) \mu_{kji} n_k n_j \\
	& \Lambda = ( \sigma_j - Q_{kj,k} ) n_j - \mathcal{D}_j ( Q_{kj} n_k ) + ( \mathcal{D}_l n_l ) Q_{kj} n_k n_j,
\end{alignat*}
and
\[
	R_i = \mu_{kji} n_k n_j, \qquad H = Q_{kj} n_k n_j,
\]
where \(\mathcal{D} \equiv n_i \partial / \partial x_i\) is the normal derivative and
\[
	\mathcal{D}_i \equiv ( \delta_{ij} - n_i n_j ) \frac{\partial}{\partial x_j}
\]
is the surface gradient.

Now, we denote with 
\[
    \mathcal{U} = ( u_i, \theta, \varphi)
\] 
the solutions of the mixed initial-boundary value problem \(\Pi\) defined by eqs.~\eqref{eq:motion-final},~\eqref{eq:dot-energy-linear},~\eqref{eq:geometric},~\eqref{eq:electric}\(_2\),~\eqref{eq:cost-eq} and the following initial conditions
\begin{alignat*}{2}
	u_i(0) &= u_i^0, \quad & \dot{u}_i(0) &= v_i^0, \\
	\theta(0) &= \theta^0, \quad & \eta(0) &= \eta^0,
\end{alignat*}
in \(\bar{B}\) and the following boundary conditions
\begin{alignat*}{2}
	u_i &= \hat{u}_i \text{ on } S_1 \times I, \quad & P_i &= \hat{P}_i \text{ on } \Sigma_1\times I, \\ 
	\mathcal{D} u_i &= \hat{d}_i \text{ on } S_2 \times I, \quad & R_i &= \hat{R}_i \text{ on } \Sigma_2\times I, \\
	\theta &= \hat{\theta} \text{ on } S_3\times I, \quad & q_i n_i &= \hat{q} \text{ on } \Sigma_3\times I, \\ 
	\varphi &= \hat{\varphi} \text{ on } S_4\times I, \quad & \Lambda &= \hat{\Lambda} \text{ on } \Sigma_4\times I, \\
	\mathcal{D}\varphi &= \hat{\xi} \text{ on } S_5\times I, \quad & H &= \hat{H} \text{ on } \Sigma_5\times I,
\end{alignat*}
with \(u_i^0\), \(v_i^0\), \(\theta^0\), \(\eta^0\), \(\hat{u}_i\), \(\hat{d}_i\), \(\hat{\theta}\), \(\hat{\varphi}\), \(\hat{\xi}\), \(\hat{P}_i\), \(\hat{R}_i\), \(\hat{q}\), \(\hat{\Lambda}\) and \(\hat{H}\) are prescribed functions and the surfaces \(S_i\) and \(\Sigma_i\) are such that
\[
	\bar{S}_i \cup \Sigma_i = \partial B \quad\quad S_i \cap \Sigma_i = \varnothing, \qquad i=1, \ldots, 5
\]
where the closure is relative to \(\partial B\).
The (external) data of the mixed initial-boundary value problem in concern are
\[
	\Gamma = \left\{ f_i, g, r, u_i^0, v_i^0, \theta^0, \eta^0, \hat{u}_i, \hat{d}_i, \hat{\theta}, \hat{\varphi}, \hat{\xi}, \hat{P}_i, \hat{R}_i, \hat{q}, \hat{\Lambda}, \hat{H} \right\}.
\]

\section{A uniqueness result}

In this section we establish a uniqueness result for a initial-boundary value problem \(\Pi \). 
To this aim, using the constitutive equation~\eqref{eq:cost-eq} we prove that
\begin{equation}
	\begin{aligned}
		 & \tau_{ij}\dot{e}_{ij}+\mu_{ijk}\dot{\kappa}_{ijk}+\dot{\sigma}_i E_i+\dot{Q}_{ji}V_{ij}+\rho\dot{\eta}(\theta+\beta\dot{\theta})= \\
		 & \qquad=\dot{W}+\dot{F}-\frac{1}{2}a^{(44)}\frac{d}{dt}(\theta+\beta\dot{\theta})^2                                                \\
		 & \qquad-\frac{1}{\beta}\left(\gamma\ddot{\theta}+a_i^{(56)}\dot{\theta}_{,i}\right)(\theta+\beta\dot{\theta}),
	\end{aligned}
	\label{eq:dotW-F-1}
\end{equation}
where \(W\) is the following quadratic form in the strain measures \(e_{ij}\) and \(\kappa_{ijk}\)
\[
	W=\frac{1}{2}a_{ijkl}^{(11)}e_{ij}e_{kl}+\frac{1}{2}a_{ijklhm}^{(22)}\kappa_{ijk}\kappa_{lhm}+a_{ijklh}^{(12)}e_{ij}\kappa_{klh}
\]
and \(F\) is a quadratic form in the variables \(E_i\), \(V_{ij}\), \(\theta+\beta\dot\theta\)
\[
	\begin{aligned}
		F & =-\frac{1}{2}a_{ij}^{(33)}E_i E_j-\frac{1}{2}a_{ijkl}^{(77)}V_{ij}V_{kl}-a_{ijk}^{(37)}E_i V_{jk} \\
		  & \quad{}-a_i^{(34)}E_i(\theta+\beta\dot{\theta})-a_{ij}^{(47)}V_{ij}(\theta+\beta\dot{\theta}).
	\end{aligned}
\]

On the other hand, taking into account eqs.~\eqref{eq:electric}\(_2\),~\eqref{eq:tau}, \eqref{eq:motion-final},~\eqref{eq:cost-eq} and~\eqref{eq:dot-energy-linear}, we have
\begin{equation}
	\begin{aligned}
		 & \tau_{ij}\dot{e}_{ij}+\mu_{ijk}\dot{\kappa}_{ijk}+\dot{\sigma}_i E_i+\dot{Q}_{ji}V_{ij}+\rho\dot{\eta}(\theta+\beta\dot{\theta}) =                 \\
		 & \quad=\left[t_{ki}\dot{u}_i+\mu_{jki}\dot{u}_{i,j}-\dot{D}_k\varphi-\dot{Q}_{kj}\varphi_{,j}-\frac{q_k}{T_0}(\theta+\beta\dot{\theta})\right]_{,k} \\
		 & \qquad{}+\rho f_i\dot{u}_i-\dot{g}\varphi+\frac{\rho r}{T_0}(\theta+\beta\dot{\theta})-\rho\ddot{u}_i\dot{u}_i                                     \\
		 & \qquad{}-\frac{1}{\beta}(\theta_{,k}+\beta\dot{\theta}_{,k})\left(a_k^{(56)}\dot{\theta}+a_{kj}^{(66)}\theta_{,j}\right)
	\end{aligned}
	\label{eq:dot-energy-3-2-1}
\end{equation}

Eqs.~\eqref{eq:dotW-F-1} and~\eqref{eq:dot-energy-3-2-1} imply
\begin{equation}
	\begin{aligned}
		 & \frac{d}{dt}\left[W+F+\frac{1}{2}\rho\dot{u}_i\dot{u}_i-\frac{1}{2}a^{(44)}(\theta+\beta\dot{\theta})^2\right]                                        \\
		 & \qquad{}-\left[t_{ki}\dot{u}_i+\mu_{jki}\dot{u}_{i,j}-\dot{D}_k\varphi-\dot{Q}_{kj}\varphi_{,j}-\frac{q_k}{T_0}(\theta+\beta\dot{\theta})\right]_{,k} \\
		 & \qquad{}-\left[\rho f_i\dot{u}_i-\dot{g}\varphi+\frac{\rho r}{T_0}(\theta+\beta\dot{\theta})\right] =                                                 \\
		 & \quad=-\frac{1}{\beta}(\theta_{,k}+\beta\dot{\theta}_{,k})\left(a_k^{(56)}\dot{\theta}+a_{kj}^{(66)}\theta_{,j}\right)                                \\
		 & \qquad+\frac{1}{\beta}\left(\gamma\ddot{\theta}+a_i^{(56)}\dot{\theta}_{,i}\right)(\theta+\beta\dot{\theta})
	\end{aligned}
	\label{eq:dot-energy-3-2-1-1}
\end{equation}

It is easy to see that
\begin{equation}
	\begin{aligned}
		 & \frac{1}{\beta}\left(\gamma\ddot{\theta}+a_i^{(56)}\dot{\theta}_{,i}\right)(\theta+\beta\dot{\theta})                                                                                                         \\
		 & \quad-\frac{1}{\beta}(\theta_{,k}+\beta\dot{\theta}_{,k})\left(a_k^{(56)}\dot{\theta}+a_{kj}^{(66)}\theta_{,j}\right)                                                                                        \\
		 & \quad=-\mathcal{P}(\dot{\theta,}\theta_{,i})+\frac{1}{2}\frac{\gamma}{\beta^2}\frac{d}{dt}(\theta+\beta\dot{\theta})^2-\frac{1}{2}\beta\frac{d}{dt}\mathcal{P}\left(-\frac{\theta}{\beta},\theta_{,i}\right)
	\end{aligned}
	\label{eq:form-p}
\end{equation}
with \(\mathcal{P}\) defined by~\eqref{eq:quadratic-form-p} and satisfying~\eqref{eq:semidef}.

Taking into account eqs.~\eqref{eq:dot-energy-3-2-1-1},~\eqref{eq:form-p} and introducing the following quadratic form in the variable \(E_i\), \(V_{ij}\), \(\theta+\beta \dot\theta\)
\[
	G = F - \frac{1}{2}\left(a^{(44)} + \frac{\gamma}{\beta^2}\right)(\theta + \beta\dot{\theta})^2
\]
we can write
\begin{equation}
	\begin{aligned}
		 & \frac{d}{dt}\left[\mathcal{E}+\frac{1}{2}\beta\mathcal{P}\left(-\frac{\theta}{\beta},\theta_{,i}\right)\right]+{}                                           \\
		 & \qquad{}-\left[t_{ki}\dot{u}_i+\mu_{jki}\dot{u}_{i,j}-\dot{D}_k\varphi-\dot{Q}_{kj}\varphi_{,j}-\frac{q_k}{T_0}(\theta+\beta\dot{\theta})\right]_{,k} \\
		 & \qquad{}-\left[\rho f_i\dot{u}_i-\dot{g}\varphi+\frac{\rho r}{T_0}(\theta+\beta\dot{\theta})\right]                                                         \\
		 & \quad={}-\mathcal{P}(\dot{\theta,}\theta_{,i})\leq0
	\end{aligned}
	\label{eq:unicita}
\end{equation}
where
\[
	\mathcal{E} = W + G + \frac{1}{2}\rho\dot{u}_i\dot{u}_i.
\]

It easily follows the next theorem

\begin{theorem}[Uniqueness]
	Assume that
	\begin{enumerate}[label=\roman*)] 
		\item \(\rho\), \(\beta\), \(T_0>0\),
		\item the constitutive coefficients satisfy the relations~\eqref{eq:symmetries},
		\item \(W\) is a positive semi-definite quadratic form,
		\item \(G\) is a positive definite quadratic form.
	\end{enumerate}
	Then, if \(S_4\) is nonempty, the initial-boundary values problem \(\Pi\) has at most one solution.
\end{theorem}

\begin{proof}
	Suppose that we have two solutions of the problem \(\Pi\).
	Then their difference \(\left(u_i, \theta, \varphi\right)\) corresponds to null data.
	By integrating eq.~\eqref{eq:unicita}, we have
	\[
		\dfrac{d}{dt}\int_B\left[\mathcal{E}+\dfrac{1}{2}\beta\mathcal{P}\left(-\dfrac{\theta}{\beta},\theta_{,i}\right)\right] \leq 0.
	\]
	The last integral must be a decreasing function with respect to time, but since it is not negative and initially null, it can be deduced that
	\begin{equation}
		\dot{u}_i = 0, \qquad \theta+\beta \dot\theta = 0, \qquad E_i = 0 \qquad \text{on} \; B \times I.
		\label{eq:funzioni nulle}
	\end{equation}
	Taking into account that eqs.~\eqref{eq:funzioni nulle}\(_{1,2}\) are linear homogeneous equations with null initial conditions and using eq.~\eqref{eq:electric}\(_2\), we conclude that
	\[
		u_i = 0, \qquad \theta=0, \qquad \varphi=\text{const} \qquad \text{on} \; B \times I.
	\]
	If \(S_4 \neq \varnothing\) we obtain the uniqueness result, in fact it is
	\[
		\varphi=0 \qquad \text{on } S_4 \times I \quad \implies \quad \varphi=0 \qquad \text{on} \; B \times I.
	\]
	\qed{}
\end{proof}

\section{Isotropic thermoelastic materials}

For the class of isotropic materials with a center of symmetry, the constitutive coefficients become
\[
	\begin{aligned}
		 & a_{ijkmnr}^{(22)} =                                                                                                       \\
		 & = \gamma_1 (
		\delta_{ij} \delta_{km} \delta_{nr} +
		\delta_{ij} \delta_{kn} \delta_{mr} +
		\delta_{ik} \delta_{jr} \delta_{mn} +
		\delta_{ir} \delta_{jk} \delta_{mn} )                                                                                        \\
		 & \quad + \gamma_{2} (
		\delta_{ik} \delta_{jm} \delta_{nr} +
		\delta_{ik} \delta_{jn} \delta_{mr} +
		\delta_{im} \delta_{jk} \delta_{nr} +
		\delta_{in} \delta_{jk} \delta_{mr} )                                                                                        \\
		 & \quad +
		\gamma_{3} \delta_{ij} \delta_{kr} \delta_{mn} +
		\gamma_{4} (
		\delta_{im} \delta_{jn} \delta_{kr} +
		\delta_{in} \delta_{jm} \delta_{kr} )                                                                                        \\
		 & \quad + \gamma_{5} (
		\delta_{im} \delta_{jr} \delta_{kn} +
		\delta_{in} \delta_{jr} \delta_{km} +
		\delta_{ir} \delta_{jm} \delta_{kn} +
		\delta_{ir} \delta_{jn} \delta_{km} ),                                                                                       \\
		 & a^{(11)}_{ijkl} = \lambda \delta_{ij} \delta_{kl} + \mu(\delta_{ik}\delta_{jl} +\delta_{il}\delta_{jk}),                  \\
		 & a^{(17)}_{ijkl} = \lambda^* \delta_{ij}\delta_{kl} + \mu^*( \delta_{ik}\delta_{jl} + \delta_{il}\delta_{jk}),             \\
		 & a^{(23)}_{ijkl} = \alpha^0\delta_{ij}\delta_{kl} + \beta^0( \delta_{ik}\delta_{jl} + \delta_{il}\delta_{jk}),             \\
		 & a^{(77)}_{ijkl} = \tilde{\lambda} \delta_{ij}\delta_{kl} + \tilde{\mu}( \delta_{ik}\delta_{jl} + \delta_{il}\delta_{jk}), \\
		 & a^{(14)}_{ij} = \alpha^{(14)}\delta_{ij}, \qquad a^{(33)}_{ij} =  \alpha^{(33)}\delta_{ij},                               \\
		 & a^{(47)}_{ij} = \alpha^{(47)}\delta_{ij}, \qquad a^{(66)}_{ij} =  \alpha^{(66)}\delta_{ij},
	\end{aligned}
\]
and
\begin{alignat*}{4}
	 & a_{ijklh}^{(12)} =0, & \quad & a_{ijk}^{(13)}=0, & \quad & a_{ijk}^{(24)}=0, \\
	 & a_{ijklh}^{(27)} =0, & \quad & a_{i}^{(34)}=0,   & \quad & a_{ijk}^{(37)} =0, & \quad & a_{i}^{(56)}=0.
\end{alignat*}
Consequently, the constitutive equations~\eqref{eq:cost-eq} reduce to
\[
	\begin{aligned}
		\tau_{ij}        & = \lambda e_{kk}\delta_{ij}+2\mu e_{ij}+\alpha^{(14)}\delta_{ij}(\theta+\beta\dot{\theta})                             \\
		                 & +\lambda^*V_{kk}\delta_{ij}+2\mu^*V_{ij},                                                                              \\
		\mu_{ijk}        & = \gamma_1(\kappa_{hhi}\delta_{jk}+2\kappa_{khh}\delta_{ij}+\kappa_{hhj}\delta_{ik})+2\gamma_2(\kappa_{ihh}\delta_{jk} \\
		                 & +\kappa_{jhh}\delta_{ik})+\gamma_3\kappa_{hhk}\delta_{ij}+2\gamma_4\kappa_{ijk}+2\gamma_5(\kappa_{kji}+\kappa_{kij})   \\
		                 & +\alpha_0\delta_{ij}E_k+\beta_0(\delta_{ik}E_j+\delta_{jk}E_i),                                                        \\
		-\sigma_i        & = \alpha_0\kappa_{kki}+2\beta_0\kappa_{ikk}+\alpha^{(33)}E_i,                                                          \\
		-\rho\eta        & = \alpha^{(14)}e_{kk}+a^{(44)}(\theta+\beta\dot{\theta})+\alpha^{(47)}V_{kk}+                                          \\
		                 & +\frac{1}{\beta}\gamma\dot{\theta}+\alpha^{(4)}                                                                        \\
		-Q_{ij}          & = \lambda^*e_{kk}\delta_{ij}+2\mu^*e_{ij}+\alpha^{(47)}\delta_{ij}(\theta+\beta\dot{\theta})                           \\
		                 & +\tilde{\lambda}V_{kk}\delta_{ij}+2\tilde{\mu}V_{ij},                                                                  \\
		-\frac{q_i}{T_0} & = \frac{1}{\beta}\alpha^{(66)}\theta_{,i}.
	\end{aligned}
\]


Using these equations and~\eqref{eq:electric}\(_2\), eqs.~\eqref{eq:motion-final} and~\eqref{eq:dot-energy-linear} can be written as follows
\begin{align*}
	 & \mu u_{i,jj}+(\lambda+\mu)u_{j,ji}+\alpha^{(14)}(\theta_{,i}+\beta\dot{\theta}_{,i})-(\lambda^*+2\mu^*)\varphi_{,ijj} \\
	 & \quad-4(\gamma_1+\gamma_2+\gamma_5)u_{j,jikk}-(\gamma_3+2\gamma_4)u_{i,jjkk}                                          \\
	 & \quad-(\alpha_0+2\beta_0)\varphi_{,ijj}+\rho f_i=\rho\ddot{u}_i,                                                      \\
	 & (\lambda^*+2\mu^*-\alpha_0-2\beta_0)u_{j,jkk}+\alpha^{(33)}\varphi_{,jj}                                              \\
	 & \quad+\alpha^{(47)}(\theta_{,jj}+\beta\dot{\theta}_{,jj})-(\tilde{\lambda}+2\tilde{\mu})\varphi_{,jjkk}=g,            \\
	 & \alpha^{(14)}\dot{u}_{j,j}+a^{(44)}\dot{\theta}+\big(a^{(44)}+\frac{1}{\beta^2}\gamma\big)\beta\ddot{\theta}          \\
	 & \quad+\frac{1}{\beta}\alpha^{(66)}\theta_{,jj}+\frac{1}{T_0}\rho r = 0.
\end{align*}
Now, we remark that the condition~\eqref{eq:semidef} of positive semi-definiteness of quadratic form \(\mathcal{P}\) is equivalent to
\[
	\frac{\gamma}{\beta} \geq 0, \qquad \frac{\alpha^{(66)}}{\beta} \geq 0.
\]
In particular, when \(\beta > 0\) eq.~\eqref{eq:semidef} is equivalent to \(\gamma \geq 0\), \(\alpha^{(66)} \geq 0\).

In the isotropic case, the quadratic form \(W\) can be expressed as the sum of two independent quadratic forms, the first one \(W_1\) in the variables \(e_{ij}\) and the second one \(W_2\) in the variables \(\kappa_{ijk}\).

The positive semi-definiteness of \(W_1\) leads to the well known conditions
\begin{equation}
	\mu\geq 0, \qquad 3\lambda + 2\mu\geq 0.
	\label{W2positivesemi-definiteness1}
\end{equation}

Regarding \(W_2\), if we choose the following ordering of the variables:
\begin{alignat*}{1}
\{ & \kappa_{221}, \kappa_{331}, \kappa_{111}, \kappa_{122}, \kappa_{133}, \\
   & \kappa_{332}, \kappa_{112}, \kappa_{222}, \kappa_{233}, \kappa_{211}, \\
   & \kappa_{113}, \kappa_{223}, \kappa_{333}, \kappa_{311}, \kappa_{322}, \\
   & \kappa_{123}, \kappa_{231}, \kappa_{312} \}
\end{alignat*}
the corresponding matrix is of type
\[
    \begin{pmatrix}
    \mathbf{A}_5 & \mathbf{0} & \mathbf{0} & \mathbf{0} \\
    \mathbf{0} & \mathbf{A}_5 & \mathbf{0} & \mathbf{0} \\
    \mathbf{0} & \mathbf{0} & \mathbf{A}_5 & \mathbf{0} \\
    \mathbf{0} & \mathbf{0} & \mathbf{0} & \mathbf{A}_3
    \end{pmatrix}
\]
where \(\mathbf{A}_5\) is a \(5\times 5\) matrix and \(\mathbf{A}_3\) is a \(3\times 3\) matrix. The eigenvalues of \(\mathbf{A}_3\) are
\begin{alignat*}{2}
    & 2(\gamma_4 - \gamma_5), \qquad && \text{multiplicity } 2 \\
    & 2(\gamma_4 + 2\gamma_5), \qquad && \text{multiplicity } 1
\end{alignat*}
The matrix \(\mathbf{A}_5\) is similar to
\[
\begin{pmatrix}
    \gamma_3 + \gamma_4 & \frac{1}{2}(2\gamma_1 + \gamma_3) & 2 (\gamma_1 + \gamma_5) & 0 & 0 \\
    2 \gamma_1 + \gamma_3 & \xi  & 2 (\gamma_1 + 2 \gamma_2) & 0 & 0 \\
    2 (\gamma_1 + \gamma_5) & \gamma_1 + 2 \gamma_2 & 2 (2 \gamma_2 + \gamma_4 + \gamma_5) & 0 & 0 \\
    0 & 0 & 0 & \xi_1 & 0 \\
    0 & 0 & 0 & 0 & \xi_2
\end{pmatrix}
\]
where
\[
    \xi = 2 (\gamma_1 + \gamma_2 + \gamma_5) + \frac{1}{2}(\gamma_3 + 2\gamma_4)
\]
and
\[
    \xi_{1,2} = \frac{1}{2} \left(3\gamma_4 + 2\gamma_5 \pm \sqrt{(\gamma_4+2\gamma_5)^2 + 16\gamma_5^2}\right)
\]

Consequently, the conditions of positive semi-definiteness of quadratic form \(W_2\) are
\begin{equation}
	\begin{aligned}
		 & \begin{aligned}
			 & \gamma_3 + \gamma_4 \geq 0,             &  & \quad \gamma_4 - \gamma_5 \geq 0, \qquad \gamma_4 + 2\gamma_5 \geq 0,    \\
			 & 2\gamma_2 + \gamma_4 + \gamma_5 \geq 0, &  & \quad 4(\gamma_1 + \gamma_2 + \gamma_5) + (\gamma_3 + 2\gamma_4) \geq 0,
		\end{aligned}                                                                                     \\
		 & (\gamma_3 + \gamma_4) (4\gamma_2 + 3\gamma_4 + 4\gamma_5) \geq (2\gamma_1 - \gamma_4)^2,                       \\
		 & (\gamma_3 + \gamma_4) (2\gamma_2 + \gamma_4 + \gamma_5) \geq 2 (\gamma_1 + \gamma_5)^2,                        \\
		 & (2\gamma_2 + \gamma_4 + \gamma_5) (\gamma_3 + 4\gamma_4 + 6\gamma_5) \geq 2(\gamma_1 - \gamma_4 - \gamma_5)^2, \\
		 & (\gamma_4 + 2\gamma_5) [(5\gamma_3 + 4\gamma_4 - 2\gamma_5) (10\gamma_2 + 3\gamma_4 + \gamma_5)                \\
		 & \qquad - 2(5\gamma_1 - \gamma_4 + 3\gamma_5)^2] \geq 0.
	\end{aligned}
	\label{W2positivesemi-definiteness2}
\end{equation}
On the other hand, the quadratic form \(G\) is positive definite if and only if
\begin{equation}
	\begin{aligned}
		 & \begin{aligned}
			 & \tilde{\mu} < 0, \qquad   &  & 3\tilde{\lambda} + 2\tilde{\mu} < 0,   \\
			 & \alpha^{(33)} > 0, \qquad &  & a^{(44)} + \frac{\gamma}{\beta^2} < 0,
		\end{aligned}                                                                           \\
		 & (3\tilde{\lambda} + 2\tilde{\mu})\left(a^{(44)} + \frac{\gamma}{\beta^2}\right)>3 (\alpha^{(47)})^2.
	\end{aligned}
	\label{Gpositivesemi-definiteness}
\end{equation}

For isotropic materials we obtain the following uniqueness result for the mixed initial-boundary values problem

\begin{theorem}
	Let us assume that
	\begin{enumerate}[label=\roman*)] 
		\item \(\rho>0\), \(T_0>0\), \(\beta>0\),  \(\gamma \geq 0\), \(a^{(66)} \geq 0\),
		\item the inequalities~\eqref{W2positivesemi-definiteness1},~\eqref{W2positivesemi-definiteness2} and~\eqref{Gpositivesemi-definiteness} hold.
	\end{enumerate}
	If \(S_4\) is nonempty, the initial-boundary values problem  considered has at most one solution.
\end{theorem}

\section{Conclusions}
We derived a theory of thermopiezoelectricity of a body in which the second gradient of displacement and the second gradient of electric potential are included in the set of independent constitutive variable.
We obtained appropriate thermodynamic restrictions and constitutive equations, with the help of an entropy production inequality proposed by \citet{Green1972a}.

For both anisotropic and isotropic materials we established the basic equations of the linear theory and obtained a uniqueness result for the mixed initial-boundary values problem.

\end{document}